\documentclass{article}
\usepackage[latin1]{inputenc}
\usepackage{amssymb}
\usepackage{amsmath,amsfonts,amsthm}
\usepackage{color}
\usepackage[]{graphicx}
\newtheorem{theorem}{Theorem}[section]
\newtheorem{lemma}[theorem]{Lemma}
\newtheorem{proposition}[theorem]{Proposition}
\newtheorem{corollary}[theorem]{Corollary}

\newtheorem{definition}[theorem]{Definition}

\newcommand{\dis}{\displaystyle}

\newcommand{\real}{\mathbb{R}}

\date{}
\begin{document}
\author{Giovany M. Figueiredo \thanks{Supported by PROCAD/CASADINHO: 552101/2011-7, CNPq/PQ  301242/2011-9 and CNPQ/CSF  200237/2012-8 }\\
\noindent Universidade Federal do Par\'a \\
\noindent Faculdade de Matem\'atica \\
\noindent CEP: 66075-110 Bel\'em - Pa , Brazil. \\
\noindent e-mail: {\tt{giovany@ufpa.br}} \\
and \\
Marcos T. O. Pimenta\thanks{Supported by FAPESP: 2012/20160-0}\, \thanks{Corresponding author}\\
\noindent Faculdade de Ci\^encias e Tecnlogia\\
\noindent Universidade Estadual Paulista - Unesp\\
\noindent 19060-900, Presidente Prudente - SP, Brazil.\\
\noindent e-mail: {\tt{pimenta@fct.unesp.br}}\\
}
\date{}

\title{Existence and multiplicity of solutions for a prescribed mean-curvature problem with critical growth} \maketitle{}

\numberwithin{equation}{section} \maketitle{}
\begin{abstract}
In this work we study an existence and multiplicity result for the following prescribed mean-curvature problem with critical growth
$$
\left\{
\begin{array}{rl}
 -\mbox{div}\biggl(\frac{\nabla u}{\sqrt{1+|\nabla u|^{2}}}\biggl) =
\lambda |u|^{q-2}u+ |u|^{2^*-2}u & \mbox{in $\Omega$}\\
u = 0 & \mbox{on $\partial \Omega$},
\end{array} \right.
$$
where $\Omega$ is a bounded
smooth domain of $\mathbb{R}^{N}$, $N\geq 3$ and $1 < q<2$.  In order to employ variational arguments, we consider an auxiliary problem which is proved to have infinitely many solutions by genus theory. A clever estimate in the gradient of the solutions of the modified problem is necessary to recover solutions of the original one.

\vspace*{0.5cm}

\noindent {\bf 2010 Mathematics Subject Classification} : 35J93, 35J62, 35J20.  \\
\noindent {\bf Key words}: prescribed mean-curvature problem,
critical exponent, variational methods.
\end{abstract}

\section{Introduction}

\hspace{0.7cm} In this work we deal with questions of existence  and multiplicity
of solutions for quasilinear problems with nonlinearity of Br\'ezis-Nirenberg type (see \cite{BrezNiremb})

$$
\left\{
\begin{array}{rl}
 -\mbox{div}\biggl(\frac{\nabla u}{\sqrt{1+|\nabla u|^{2}}}\biggl) =
\lambda |u|^{q-2}u+ |u|^{2^*-2}u & \mbox{in $\Omega$}\\
u = 0 & \mbox{on $\partial \Omega$},
\end{array} \right.
 \eqno{(P_{\lambda})}
$$
\noindent where $\Omega\subset\mathbb{R}^{N}$ is a bounded smooth
domain, $\lambda > 0$, $1 < q<2$ and  $2^{*}=\frac{2N}{N-2}$. This kind of problem has applications not just to describe a surface given by  $u(x)$, whose mean curvature is described by the right hand side of $(P_\lambda)$, but also in {\it capillarity theory} where when the nonlinearity is replaced by  $\kappa u$, the resultant equation describe the equilibrium of a liquid surface with constant surface tension in a uniform gravity field (see p. 262 in \cite{Trudinger}).

Problems like $(P_\lambda)$ has been intensively studied over the last decades. In the work \cite{Coffman}, the authors studied a related subcritical problem in which they obtained positive solutions. In the recent work \cite{Bonheure1}, Bonheure, Derlet and Valeriola have studied a purely subcritical version of $(P_\lambda)$, where they proved the existence and multiplicity of nodal $H^1_0(\Omega)$ solutions, to sufficiently large values of $\lambda$. They overcame the difficulty in working in the $BV(\Omega)$ space, which is the natural functional space to treat $(P_\lambda)$, by doing a truncation in the degenerate part of the mean-curvature operator in order to make possible construct a variational framework in the Sobolev space $H^1_0(\Omega)$. Nevertheless, this truncation requires sharp estimates on the gradient of the solutions, in order to prove that the solutions of the modified problem in fact are solutions of the original one.

When $\Omega = \mathbb{R}^N$ and the nonlinearity is substituted by $u^q$, i.e., the Gidas-Spruck analogue for the mean-curvature operator, Ni and Serrin in \cite{NiSerrin1,NiSerrin2} has proved that if $1 < q < \frac{N}{N-2}$ no positive solution exist, while for $q \geq 2^*-1$ there exist infinitely many solutions. In the range $\frac{N}{N-2} < q < 2^*-1$ some contributions  has been given by Cl\'ement et al in \cite{Clement} and by Del Pino and Guerra in \cite{DelPino}, where in the latter the authors prove that many positive solutions do exist if $q < 2^*-1$ is sufficiently close to $2^*-1$.

Still in the case $\Omega = \mathbb{R}^N$ but with nonlinearity given by $\lambda u + u^p$, Peletier and Serrin in \cite{Peletier} succeed in proving the existence of positive radial solutions when $\lambda < 0$ is small enough and $p$ is subcritical. In the case $\lambda > 0$, they stated there is no regular solution to that problem no matter how much small or large $p$ is.

In this work, because of the boundedness of $\Omega$, we prove a result in a strike opposition of that in \cite{Peletier}, in which we obtain the existence of infinitely many regular solutions of $(P_\lambda)$, for small enough $\lambda > 0$. More specifically, we prove the following result.

\begin{theorem}\label{Teorema1}
If $1<q<2$, then there exists $\lambda^{*} > 0$ such
that if $0 < \lambda < \lambda^*$, $(P_{\lambda})$ has infinitely many solutions. Moreover, if $u_{\lambda}$ is a
solution of $(P_{\lambda})$, then $u_{\lambda} \in
H^1_0(\Omega) \cap C^{1,\alpha}(\overline{\Omega})$ with $\alpha
\in (0,1)$, and
$$
\displaystyle\lim_{\lambda\rightarrow 0}\|u_{\lambda}\| = \displaystyle\lim_{\lambda\rightarrow 0}\|u_{\lambda}\|_{\infty} = \displaystyle\lim_{\lambda\rightarrow 0}\|\nabla
u_{\lambda}\|_{\infty}=0,
$$
where $\|\cdot\|$ is the Sobolev norm in $H^1_0(\Omega)$.
\end{theorem}

Our approach follows the main ideas of Bonheure {\it et al} in \cite{Bonheure1}, in order to make possible consider a related modified problem in $H^1_0(\Omega)$. Afterwards, to get solutions of the modified problem we apply {\it Krasnoselskii genus theory} in the same way that Azorero and Alonso in \cite{GP}. Finally, we use the Moser iteration technique and a regularity result by Lewy and Stampacchia in \cite{Stampacchia} to get decay in $\lambda$ of the gradient of the solutions, which will imply that the solutions of the modified problem in fact are solutions of the original one.

The paper is organized as follows. In the second section we present the auxiliary problem and the variational framework. In the third one we make a brief review of Genus theory. In the fourth we prove some technical results which imply on the existence of infinitely many solutions to the auxiliary problem. The last one is dedicated to present the proof of the main result, which consists in estimates in $L^\infty(\Omega)$ norm of the gradient of solutions.

\section{The auxiliary problem and variational framework}

\hspace{0.7cm} Let us consider $r \geq
0$ and $\delta >0$ and a function $ \eta \in C^{1}([r,
r+\delta])$ such that
$$
\eta (r)=\frac{1}{\sqrt{1+r}}, \ \ \eta
(r+\delta)=\frac{1}{\sqrt{1+r+\delta}},
$$
$$
\eta' (r)=-\frac{1}{2\sqrt{(1+r)^{3}}} \ \ \mbox{and} \ \ \eta' (r+
\delta)=0.
$$

Now we define

$$
a(t):= \left\{
\begin{array}{cl}
\frac{1}{\sqrt{1+t}}, & \mbox{if } \  0 \leq t \leq r, \\
\eta(t) & \mbox{if }  \ r\leq t \leq r + \delta, \\
K_{0}=\frac{1}{\sqrt{1+r+\delta }},&  \mbox{if } \   t \geq
r+\delta.
\end{array}
\right.
$$

Note that $ a \in  C^{1}([0, \infty))$ is decreasing  and $K_{0}
\leq a(t)\leq 1$ for $t \in [0, \infty)$. Let us fix $r>0$
such that
\begin{eqnarray}\label{pradacaerto1}
\frac{2}{2^{*}}< K_{0}< 1.
\end{eqnarray}

The proof of the Theorem \ref{Teorema1} is
based on a careful study of solutions of the following auxiliary
problem
$$
\left\{
\begin{array}{rl}
 -\mbox{div}(a(|\nabla u|^2)\nabla u) =
\lambda |u|^{q-2}u+ |u|^{2^*-2}u & \mbox{in $\Omega$}\\
u = 0 & \mbox{on $\partial \Omega$},
\end{array} \right.
\leqno{(T_{\lambda})}
$$

We say that $u \in H^{1}_{0}(\Omega)$ is a weak solution $(T_{\lambda})$ if it verifies
\begin{eqnarray*}
\dis\int_{\Omega}a(|\nabla u|^{2})\nabla u \nabla \phi \ dx =
\lambda\dis\int_{\Omega}|u|^{q-2}u\phi
 \ dx +\dis\int_{\Omega}|u|^{2^{*}-2}u\phi  \ dx,
\end{eqnarray*}
for all $\phi \in  H^{1}_{0}(\Omega)$. Let us consider $H^1_0(\Omega)$ with its usual norm $\|u\|= \left(\int_\Omega|\nabla u|^2\right)^\frac{1}{2}$ and define the
$C^1$-functional $I_{\lambda}:  H^{1}_{0}(\Omega)\to \mathbb{R}$ by
$$
I_{\lambda}(u) = \dis\frac{1}{2}\dis\int_{\Omega}A(|\nabla u|^{2}) \
dx - \dis\frac{\lambda}{q}\dis\int_{\Omega}|u|^{q} \ dx
-\dis\frac{1}{2^{*}}\dis\int_{\Omega}|u|^{2^{*}} \ dx ,
$$
where $A(t)=\displaystyle\int^{t}_{0}a(s) \ ds$. Note that
\begin{eqnarray*}
I_{\lambda}'(u)\phi = \displaystyle\int_{\Omega}a(|\nabla u|^{2})
\nabla u \nabla \phi \ dx -
\lambda\displaystyle\int_{\Omega}|u|^{q-2}u\phi
 \ dx-\displaystyle\int_{\Omega}|u|^{2^{*}-2}u\phi  \ dx,
\end{eqnarray*}
for all $\phi \in   H^{1}_{0}(\Omega)$ and then, critical points of
$I_{\lambda}$ are weak solutions of $(T_{\lambda})$.

In order to use variational methods, we first derive some results
related to the Palais-Smale  compactness condition.

We say that a sequence $(u_{n})\subset   H^{1}_{0}(\Omega)$ is a
$(PS)_{c_\lambda}$ sequence for $I_{\lambda}$ if
\begin{eqnarray}\label{***}
I_{\lambda}(u_{n})\rightarrow c_{\lambda} \ \mbox{and} \
\|I_{\lambda}'(u_{n})\|_{H^{-1}(\Omega)}\rightarrow 0, \quad \mbox{as $n\to \infty$}
\end{eqnarray}
where $$ c_{\lambda} = \displaystyle\inf_{\pi \in \Gamma}
\displaystyle\max_{t \in [0,1]} I_{\lambda}(\pi(t))>0$$ and
$$
\Gamma := \{ \pi\in C([0,1],H^{1}_{0}(\Omega)) : \pi(0)=0,
~I_{\lambda}(\pi(1)) < 0\}.
$$

If (\ref{***}) implies the existence of a subsequence $(u_{n_{j}})
\subset (u_{n})$ which converges in $H^1_0(\Omega)$, we say that $I_{\lambda}$
satisfies the Palais-Smale condition on the level $c_\lambda$.

\section{Genus theory}

\hspace{0.7cm} We start by considering some basic facts on the Krasnoselskii
genus theory that we will use in the proof of Theorem \ref{Teorema1}.

Let $E$ be a real Banach space. Let us denote by $\mathfrak{A}$ the
class off all closed subsets  $A\subset E\setminus \{0\}$ that are
symmetric with respect to the origin, that is, $u\in A$ implies
$-u\in A$.

\begin{definition}
Let $A\in \mathfrak{A}$. The Krasnoselskii genus $\gamma(A)$ of $A$
is defined as being the least positive integer $k$ such that there
is an odd mapping $\phi \in C(A,\real^{k})$ such that $\phi(x)\neq
0$ for all $x\in A$. When such number does not exist we set
$\gamma(A)=\infty$. Furthermore, by definition,
$\gamma(\emptyset)=0$.
\end{definition}

In the sequel we will establish only the properties of the genus
that will be used through this work. More informations on this
subject may be found in \cite{kras}.
\begin{theorem}
Let $E={\real}^{N}$ and $\partial\Omega$ be the boundary of an open,
symmetric and bounded subset $\Omega \subset {\real}^{N}$ such that $0
\in \Omega$. Then $\gamma(\partial\Omega)=N$.
\end{theorem}

\begin{corollary}\label{esfera}
$\gamma(S^{N-1})=N$.
\end{corollary}

\begin{proposition}\label{paracompletar}
If $K \in \mathfrak{A}$, $0 \notin K$ and $\gamma(K) \geq 2$, then
$K$ has infinitely many points.
\end{proposition}


\section{Technical results}


\hspace{0.7cm} The genus theory requires that the functional $I_{\lambda}$ is
bounded from below. Since this is not the case, it is necessary to work with a related functional, which will be done employing some ideas contained in \cite{GP}.

In the light of the Proposition \ref{paracompletar}, it seems to be useful prove that the set of critical points of the related functional has genus greater than 2, in order to obtain infinitely many solutions of $(T_{\lambda})$. 

Lets gonna present the way in which we truncate the function $I_\lambda$ .
From (\ref{pradacaerto1}) and Sobolev's embedding, we get
\begin{eqnarray*}
I_{\lambda}(u) \geq \frac{K_{0}}{2}\|u\|^{2} -\frac{\lambda}{q
S_{q}^{q/2}}\|u\|^{q} -
\frac{1}{2^{*}S^{2^{*}/2}}\|u\|^{2^{*}}=g(\|u\|^{2}),
\end{eqnarray*}
$S$ and $S_{q}$ are, respectively, the best constants of the Sobolev's embeddings
$H^1_0(\Omega)\hookrightarrow L^{2^*}(\Omega)$ and $H^1_0(\Omega) \hookrightarrow L^{q}(\Omega)$ and
\begin{eqnarray}\label{gabriel1}
g(t)=\frac{K_{0}}{2} t -\displaystyle\frac{\lambda}{q
S_{q}^{q/2}}t^{q/2} -
\displaystyle\frac{1}{2^{*}S^{2^{*}/2}}t^{2^{*}/2}.
\end{eqnarray}

Hence, there exists $\tau_{1}>0$ such that, if $\lambda \in (0,
\tau_{1})$, then $g$ attains its positive maximum.

Let $R_{0} < R_{1}$ the roots of $g$. We have that
$R_{0}=R_{0}(\tau_{1})$ and the following result holds:

\begin{lemma}\label{comportamentoassimtotico100}
\begin{eqnarray}\label{comportamentoassimtotico1}
R_{0}(\tau_{1})\rightarrow 0 \ \ \mbox{as} \ \ \lambda\rightarrow
0.
\end{eqnarray}
\end{lemma}
\noindent\textbf{Proof:} From $g(R_{0}(\tau_{1}))=0$ and
$g'(R_{0}(\tau_{1}))>0$, we have
\begin{equation}\label{lambda1}
\frac{K_{0}}{2}R_{0}(\tau_{1})=\displaystyle\frac{\lambda}{q
S_{q}^{q/2}}R_{0}(\tau_{1})^{q/2}+
\displaystyle\frac{1}{2^{*}S^{2^{*}/2}}R_{0}(\tau_{1})^{2^{*}/2}
\end{equation}
and
\begin{equation}\label{lambda2}
\frac{K_{0}}{2}>\displaystyle\frac{\lambda}{2q
S_{q}^{q/2}}R_{0}(\tau_{1})^{(q-2)/2}+
\displaystyle\frac{1}{2S^{2^{*}/2}}R_{0}(\tau_{1})^{(2^{*}-2)/2},
\end{equation}
for all $\lambda\in (0,\tau_{1})$. From, $(\ref{lambda1})$ we
conclude that $R_{0}(\tau_{1})$ is bounded. Suppose that
$R_{0}(\tau_{1})\to R_{0}>0$ as $\lambda\to 0$. Then,
\begin{equation}
\frac{K_{0}}{2}=
\displaystyle\frac{1}{2^{*}S^{2^{*}/2}}R_{0}(\tau_{1})^{2^{*}-2/2}
\end{equation}
and
\begin{equation}\label{lambda2}
\frac{K_{0}}{2}\geq
\displaystyle\frac{1}{2S^{2^{*}/2}}R_{0}(\tau_{1})^{2^{*}-2/2},
\end{equation}
which is a contradiction, because $2^{*}>2$. Therefore $R_{0}=0$.
\hfill\rule{2mm}{2mm}

\vspace{0.5cm}

We consider $\tau_{1}$ such that $R_{0}\leq r$ and we modify the functional $I_{\lambda}$ in the following way. Take $\phi \in C^{\infty}([0,+\infty))$, $0\leq
\phi\leq 1$  such that $\phi(t)=1$
if $t\in [0,R_{0}]$ and $\phi(t)=0$ if $t\in [R_{1}, +\infty)$. Now,
we consider the truncated functional
$$
J_{\lambda}(u) = \dis\frac{1}{2}\dis\int_{\Omega}A(|\nabla u|^{2}) \
dx - \displaystyle\frac{\lambda}{q}\displaystyle\int_{\Omega}|u|^{q}
\ dx -
\phi(\|u\|^{2})\displaystyle\frac{1}{2^{*}}\displaystyle\int_{\Omega}|u|^{2^{*}}
\ dx.
$$

Note that $J_{\lambda} \in C^{1}(H^{1}_{0}(\Omega),\mathbb{R})$ and,
as in (\ref{gabriel1}), $J_{\lambda}(u)\geq
\overline{g}(\|u\|^{2})$, where
$$
\overline{g}(t)= \frac{K_{0}}{2} t -\displaystyle\frac{\lambda}{q
S_{q}^{q/2}}t^{q/2} -
\phi(t)\displaystyle\frac{1}{2^{*}S^{2^{*}/2}}t^{2^{*}/2}.
$$

Let us remark that if $\|u\|^{2}\leq R_{0}$, then
$J_{\lambda}(u)=I_{\lambda}(u)$ and if $\|u\|^{2}\geq R_{1}$, then
$J_{\lambda}(u)= \dis\frac{1}{2}\dis\int_{\Omega}A(|\nabla u|^{2}) \
dx- \displaystyle\frac{\lambda}{q}\displaystyle\int_{\Omega}|u|^{q}
\ dx$, which implies that $J_{\lambda}$ is coercive and hence bounded from below.

Now we show that $J_{\lambda}$ satisfy the local Palais-Smale
condition. For this, we need the following technical result, which
is analogous of Lemma 4.2 in \cite{GP}.

\begin{lemma}\label{nivelbaixo}
Let $(u_{n}) \subset H^{1}_{0}(\Omega)$ be  a bounded sequence such
that
$$
I_{\lambda}(u_{n})\rightarrow c_{\lambda} \ \ \mbox{and} \ \
I_{\lambda}'(u_{n})\rightarrow 0.
$$
If
\begin{eqnarray*}
c_{\lambda}&< & (\displaystyle\frac{K_{0}}{2}-
\frac{1}{2^{*}})K_{0}^{(N-2)/2}S^{N/2} \\
&- & \lambda
(\frac{1}{q}-\frac{1}{2^{*}})|\Omega|^{\frac{(2^{*}-q)}{2}}\biggl[\frac{q}{2^{*}}\lambda\bigl(\frac{1}{q}-\frac{1}{2^{*}}\bigl)|\Omega|^{(2^{*}-q)/2^{*}}
\biggl(\bigl(\displaystyle\frac{K_{0}}{2}-
\frac{1}{2^{*}})\frac{1}{S^{2^{*}/2}}\biggl)^{-1}\biggl]^{\frac{q}{(2^{*}-q)}}
\end{eqnarray*}
hold, then up to a subsequence $(u_{n})$ is strongly
convergent in  $H^{1}_{0}(\Omega)$.
\end{lemma}
\noindent\textbf{Proof:}   Taking a subsequence, we may suppose that
\begin{eqnarray*}
|\nabla u_n|^2 \rightharpoonup |\nabla u|^2 + \sigma \ \  \text{ and
} \ \ |u_n|^{2^{*}} \rightharpoonup |u|^{2^{*}} +
\nu \quad \mbox{in the weak* sense of measures.}
\end{eqnarray*}

Using the concentration compactness-principle due to Lions (cf.
\cite[Lemma 2.1]{Lio2}), we obtain an at most countable index set
$\Lambda$, sequences $(x_i) \subset \Omega$, $(\mu_i), (\sigma_i),
(\nu_i), \subset [0,\infty)$, such that
\begin{equation}
\nu  =  \sum_{i \in \Lambda}\nu_{i}\delta_{x_{i}}, ~~~ \sigma\geq
\sum_{i \in \Lambda}\sigma_{i}\delta_{x_{i}}
 ~~\text{ and }~~S
\nu_{i}^{2/2^{*}}\leq \sigma_{i},
 \label{lema_infinito_eq11}
\end{equation}
for all $i \in\Lambda$, where $\delta_{x_i}$ is the Dirac mass at
$x_i \in \Omega$.

Now we claim that $\Lambda=\emptyset$. Arguing by contradiction,
assume that $\Lambda\neq\emptyset$ and fix $i \in \Lambda$. Consider
$\psi \in C_0^{\infty}(\Omega,[0,1])$ such that $\psi \equiv 1$ on
$B_1(0)$, $\psi \equiv 0$ on $\Omega \setminus B_2(0)$ and $|\nabla
\psi|_{\infty} \leq 2$. Defining $\psi_{\varrho}(x) :=
\psi((x-x_i)/\varrho)$ where  $\varrho>0$, we have that
$(\psi_{\varrho}u_n)$ is bounded. Thus
$I_{\lambda}'(u_n)(\psi_{\varrho}u_n) \to 0$, that is,

$$
\begin{array}{lcl}
&& \displaystyle\int_{\Omega}a(|\nabla u_{n}|^{2})u_{n} \nabla u_n
\nabla \psi_{\varrho} \ dx + \displaystyle\int_{\Omega} a(|\nabla
u_{n}|^{2})\psi_{\varrho}|\nabla u_n|^2 \ dx
\\ &=& \lambda\displaystyle\int_{\Omega}|u_n|^{q}\psi_{\varrho} \
dx+ \displaystyle\int_{\Omega} \psi_{\varrho}|u_n|^{2^{*}} \ dx +
o_n(1).
\end{array}
$$

Since $supp(\psi_{\varrho}) \subset B_{2\varrho}(x_{i})$, we
obtain
$$
\biggl|\displaystyle\int_{\Omega}u_{n}\nabla u_{n} \nabla
\psi_{\varrho} \ dx \biggl|\leq \int_{B_{2\rho}(x_{i})} |\nabla
u_{n}||u_{n} \nabla \psi_{\varrho}| \ dx.
$$

By H\"older inequality and the fact that the sequence $(u_{n})$ is
bounded in $H^{1}_{0}(\Omega)$ we have
$$
\biggl|\displaystyle\int_{\Omega}u_{n}\nabla u_{n} \nabla
\psi_{\varrho} \ dx \biggl|\leq C \left( \int_{B_{2\varrho}(x_{i})}|u_{n}
\nabla \psi_{\varrho}|^{2} \ dx\right)^{1/2}.
$$

By the Dominated Convergence Theorem
$\displaystyle\int_{B_{2\varrho}(x_{i})} |u_{n}\nabla
\psi_{\varrho}|^{2} \ dx \to 0$ as $n \to +\infty$ and $\varrho \to
0$. Thus, we obtain
\begin{equation*}
\displaystyle\lim_{\varrho\rightarrow
0}\left[\displaystyle\lim_{n\rightarrow \infty}\displaystyle\int_{\Omega}
u_{n}\nabla u_n  \nabla \psi_{\varrho} \ dx \right]=0.
\end{equation*}

Since $0<K_{0}\leq a(t) \leq 1$, for all $t \in \mathbb{R}$, we get
$$
\displaystyle\lim_{\varrho\rightarrow 0}\lim_{n\rightarrow
\infty}\left[\displaystyle \displaystyle\int_{\Omega}a(|\nabla
u_{n}|^{2})u_{n} \nabla u_n \nabla \psi_{\varrho} \ dx \right]=0.
$$

Moreover, similar arguments applies in order to obtain
$$
\displaystyle\lim_{\varrho\rightarrow 0}\lim_{n\rightarrow
\infty}\left[\displaystyle\int_{\Omega}\psi_{\varrho}|u_n|^{q} \ dx\right]=0.
$$

Thus, we have
$$
K_{0}\int_{\Omega} \psi_{\varrho}\textrm{d}\sigma \leq \int_{\Omega}
\psi_{\varrho}\textrm{d}\nu+o_{\varrho}(1).
$$
Letting $\varrho \to 0$ and using standard theory of Radon measures,
we conclude that $K_{0}\sigma_i \leq \nu_i$. It follows from
(\ref{lema_infinito_eq11}) that

\begin{eqnarray} \label{lemafinitoeq22}
\sigma_{i} \geq K_{0}^{(N-2)/2}S^{N/2}.
\end{eqnarray}

Now we shall prove that the above expression cannot occur, and
therefore the set $\Lambda$ is empty. Indeed, if for some $i \in \Lambda$ (\ref{lemafinitoeq22}) hold, then
\begin{eqnarray*}
c_{\lambda} &=& I_{\lambda}(u_n) - \displaystyle\frac{1}{2^{*}}
I_{\lambda}'(u_n)u_n + o_n(1)
\end{eqnarray*}
implies that
\begin{eqnarray*}
c_{\lambda} \geq (\displaystyle\frac{K_{0}}{2}-
\frac{1}{2^{*}})\displaystyle\int_{\Omega}|\nabla u_{n}|^{2} \ dx-
\lambda \displaystyle
(\frac{1}{q}-\frac{1}{2^{*}})\displaystyle\int_{\Omega} |u_n|^{q} \
dx.
\end{eqnarray*}

Since $\frac{2}{2^{*}}< K_{0}< 1$ ( see (\ref{pradacaerto1})),
letting $n\to\infty$ we get
\begin{eqnarray*}
c_{\lambda} \geq (\displaystyle\frac{K_{0}}{2}-
\frac{1}{2^{*}})\sigma_{i}+ (\displaystyle\frac{K_{0}}{2}-
\frac{1}{2^{*}})\displaystyle\int_{\Omega}|\nabla u|^{2} \ dx -
\lambda \displaystyle
(\frac{1}{q}-\frac{1}{2^{*}})\displaystyle\int_{\Omega} |u|^{q} \
dx.
\end{eqnarray*}

Hence,

\begin{eqnarray*}
c_{\lambda} \geq (\displaystyle\frac{K_{0}}{2}-
\frac{1}{2^{*}})K_{0}^{(N-2)/2}S^{N/2}+
(\displaystyle\frac{K_{0}}{2}-
\frac{1}{2^{*}})\displaystyle\int_{\Omega}|\nabla u|^{2} \ dx-
\lambda \displaystyle
(\frac{1}{q}-\frac{1}{2^{*}})\displaystyle\int_{\Omega} |u|^{q} \
dx.
\end{eqnarray*}

By H\"older's inequality and Sobolev's embedding we obtain
\begin{eqnarray*}
c_{\lambda}  &\geq & (\displaystyle\frac{K_{0}}{2}-
\frac{1}{2^{*}})K_{0}^{(N-2)/2}S^{N/2}+
(\displaystyle\frac{K_{0}}{2}-
\frac{1}{2^{*}})\frac{1}{S^{2^{*}/2}}\displaystyle\int_{\Omega}|
u|^{2^{*}} \ dx \\
&-& \lambda \displaystyle
(\frac{1}{q}-\frac{1}{2^{*}})|\Omega|^{\frac{(2^{*}-q)}{2^{*}}}\bigg(\displaystyle\int_{\Omega}
|u|^{2^{*}} \ dx\biggl)^{q/2^{*}}.
\end{eqnarray*}

Note that
$$
f(t)=(\displaystyle\frac{K_{0}}{2}-
\frac{1}{2^{*}})\frac{1}{S^{2^{*}/2}}t^{2^{*}}-\lambda
(\frac{1}{q}-\frac{1}{2^{*}})|\Omega|^{\frac{(2^{*}-q)}{2^{*}}}t^{q}
$$
is a continuous function that attains its absolute minimum, for
$t>0$, at the point

$$
\alpha_{0}=\biggl[\frac{q}{2^{*}}\lambda\bigl(\frac{1}{q}-\frac{1}{2^{*}}\bigl)|\Omega|^{(2^{*}-q)/2^{*}}
\biggl(\bigl(\displaystyle\frac{K_{0}}{2}-
\frac{1}{2^{*}})\frac{1}{S^{2^{*}/2}}\biggl)^{-1}\biggl]^{\frac{1}{(2^{*}-q)}}.
$$

Then,
\begin{eqnarray*}
c_{\lambda} \geq (\displaystyle\frac{K_{0}}{2}-
\frac{1}{2^{*}})K_{0}^{(N-2)/2}S^{N/2}+
(\displaystyle\frac{K_{0}}{2}-
\frac{1}{2^{*}})\frac{1}{S^{2^{*}/2}}\alpha^{2^{*}}_{0}- \lambda
(\frac{1}{q}-\frac{1}{2^{*}})|\Omega|^{\frac{(2^{*}-q)}{2}}\alpha_{0}^{q}.
\end{eqnarray*}

So
\begin{eqnarray*}
c_{\lambda} \geq (\displaystyle\frac{K_{0}}{2}-
\frac{1}{2^{*}})K_{0}^{(N-2)/2}S^{N/2} -  \lambda
(\frac{1}{q}-\frac{1}{2^{*}})|\Omega|^{\frac{(2^{*}-q)}{2}}\alpha_{0}^{q}.
\end{eqnarray*}

Thus, we conclude that
\begin{eqnarray*}
c_{\lambda} &\geq &(\displaystyle\frac{K_{0}}{2}-
\frac{1}{2^{*}})K_{0}^{(N-2)/2}S^{N/2} \\
&-&  \lambda
(\frac{1}{q}-\frac{1}{2^{*}})|\Omega|^{\frac{(2^{*}-q)}{2}}\biggl[\frac{q}{2^{*}}\lambda\bigl(\frac{1}{q}-\frac{1}{2^{*}}\bigl)|\Omega|^{(2^{*}-q)/2^{*}}
\biggl(\bigl(\displaystyle\frac{K_{0}}{2}-
\frac{1}{2^{*}})\frac{1}{S^{2^{*}/2}}\biggl)^{-1}\biggl]^{\frac{q}{(2^{*}-q)}},
\end{eqnarray*}
which is a contradiction. Thus $\Lambda$ is empty and it follows
that $u_n \to u$ in $L^{2^{*}}(\Omega)$. Thus, up to a subsequence,
$$
\|u_{n}-u\|^{2} \leq
\frac{1}{K_{0}}\displaystyle\int_{\Omega}a(|\nabla
u_{n}|^{2})|\nabla u_{n} - \nabla u|^{2} = I_{\lambda}(u_{n})u_{n}-
I_{\lambda}(u_{n})u+ o_{n}(1)=o_{n}(1).
$$
\hfill\rule{2mm}{2mm}

By the Lemma \ref{nivelbaixo}  we conclude that, there exists
$\tau_{2}>0$ such that, for all $\lambda \in (0,\tau_{2})$ we get
\begin{eqnarray*}
&&(\displaystyle\frac{K_{0}}{2}-
\frac{1}{2^{*}})K_{0}^{(N-2)/2}S^{N/2} \\
&-&  \lambda
(\frac{1}{q}-\frac{1}{2^{*}})|\Omega|^{\frac{(2^{*}-q)}{2}}\biggl[\frac{q}{2^{*}}\lambda\bigl(\frac{1}{q}-\frac{1}{2^{*}}\bigl)|\Omega|^{(2^{*}-q)/2^{*}}
\biggl(\bigl(\displaystyle\frac{K_{0}}{2}-
\frac{1}{2^{*}})\frac{1}{S^{2^{*}/2}}\biggl)^{-1}\biggl]^{\frac{q}{(2^{*}-q)}}
>0
\end{eqnarray*}
and, hence, if $(u_{n})$ is a bounded sequence such that
$I_{\lambda}(u_{n})\rightarrow c$,  $I_{\lambda}'(u_{n})\rightarrow
0$ with $c<0$, then $(u_{n})$ has a convergent subsequence.

\begin{lemma}\label{faltava}
If $J_{\lambda}(u) <0$, then $\|u\|^{2}< R_{0}\leq r$ and
$J_{\lambda}(v)=I_{\lambda}(v)$, for all $v$ in a small enough
neighborhood of $u$. Moreover, $J_{\lambda}$ verifies a local
Palais-Smale condition for $c <0$.
\end{lemma}
\noindent\textbf{Proof:} Since $\overline{g}(\|u\|^{2})\leq
J_{\lambda}(u) <0$, then $\|u\|^{2} < R_{0}\leq r$. By the choice of
$\tau_{1}$ in (\ref{comportamentoassimtotico1}) we have that
$J_{\lambda}(u)=I_{\lambda}(u)$. Moreover, since $J_{\lambda}$ is continuous, we conclude that
$J_{\lambda}(v)=I_{\lambda}(v)$, for all $v \in B_{R_{0}/2}(0)$.
Besides, if $(u_{n})$ is a sequence such that
$J_{\lambda}(u_{n})\rightarrow c<0$ and
$J_{\lambda}'(u_{n})\rightarrow 0$ as $n\to \infty$, then for $n$ sufficiently large
$I_{\lambda}(u_{n})=J_{\lambda}(u_{n})\rightarrow c<0$ and
$I_{\lambda}'(u_{n})=J_{\lambda}'(u_{n})\rightarrow 0$ as $n \to \infty$. Since
$J_{\lambda}$ is coercive, we get that $(u_{n})$ is bounded in
$H^{1}_{0}(\Omega)$. From Lemma \ref{nivelbaixo}, for all $\lambda
\in (0,\tau_{2})$, we obtain
\begin{eqnarray*}
&&c<0 < (\displaystyle\frac{K_{0}}{2}-
\frac{1}{2^{*}})K_{0}^{(N-2)/2}S^{N/2} \\
&-&  \lambda
(\frac{1}{q}-\frac{1}{2^{*}})|\Omega|^{\frac{(2^{*}-q)}{2}}\biggl[\frac{q}{2^{*}}\lambda\bigl(\frac{1}{q}-\frac{1}{2^{*}}\bigl)|\Omega|^{(2^{*}-q)/2^{*}}
\biggl(\bigl(\displaystyle\frac{K_{0}}{2}-
\frac{1}{2^{*}})\frac{1}{S^{2^{*}/2}}\biggl)^{-1}\biggl]^{\frac{q}{(2^{*}-q)}}
\end{eqnarray*}
and hence, up to a subsequence $(u_{n})$ is strongly convergent in
$H^{1}_{0}(\Omega)$. \hfill\rule{2mm}{2mm}

Now, we construct an appropriate minimax sequence of negative
critical values.

\begin{lemma}\label{minimax}
Given $k \in \mathbb{N}$, there exists $\epsilon = \epsilon(k)>0$
such that
$$
\gamma(J_{\lambda}^{-\epsilon}) \geq k,
$$
where $J_{\lambda}^{-\epsilon}=\{u \in H^{1}_{0}(\Omega):
J_{\lambda}(u) \leq -\epsilon\}$.
\end{lemma}
\noindent\textbf{Proof:} Consider $k \in \mathbb{N}$ and let $X_{k}$
be a k-dimensional subspace of $H^{1}_{0}(\Omega)$. Since in $X_k$ all norms are equivalent, there
exists $C(k)>0$ such that
$$
-C(k)\|u\|^{q}\geq - \displaystyle\int_{\Omega}|u|^{q} \ dx,
$$
for all $u \in X_{k}$.

We now use the inequality above  to conclude that
\begin{eqnarray*}
J_{\lambda}(u)\leq \frac{1}{2}\| u\|^{2}-\frac{C(k)}{q}\|u\|^{q}= \|
u \|^{q}\left(\frac{1}{2}\| u\|^{2-q}-\frac{C(k)}{q}\right).
\end{eqnarray*}

Considering $R>0$ sufficiently small, there exists
$\epsilon=\epsilon(R)>0$ such that
$$
J_{\lambda}(u)<-\epsilon < 0,
$$
for all $u\in {\mathcal{S}_R}=\{u\in X_k; \| u \|=R \}$. Since $X_k$
and $\mathbb{R}^k$ are isomorphic and $\mathcal{S}_R$ and $S^{k-1}$
are homeomorphic, we conclude from Corollary \ref{esfera} that
$\gamma(\mathcal{S}_R)=\gamma(S^{k-1})=k$. Moreover, once that
${\mathcal{S}_R} \subset J_{\lambda}^{-\epsilon}$ and
$J_{\lambda}^{-\epsilon}$ is symmetric and closed,  we have
$$
k= \gamma ({\mathcal{S}_R})\leq \gamma( J_{\lambda}^{-\epsilon}).
$$
\hfill\rule{2mm}{2mm}

We define now, for each $k \in \mathbb{N}$, the sets
$$
\Gamma_{k}=\{C \subset H: C \ \ \mbox{is closed}, C=-C \ \
\mbox{and} \ \ \gamma(C) \geq k\},
$$
$$
K_{c}=\{u \in H: J_{\lambda}'(u)=0 \ \ \mbox{and} \ \
J_{\lambda}(u)=c\}
$$
and the number
$$
c_{k}=\displaystyle\inf_{C\in \Gamma_{k}}\displaystyle\sup_{u \in
C}J_{\lambda}(u).
$$

\begin{lemma}\label{minimax1}
Given $k \in \mathbb{N}$, the number $c_{k}$ is negative.
\end{lemma}
\noindent\textbf{Proof:} From Lemma \ref{minimax}, for each $k\in
\mathbb{N}$ there exists $\epsilon >0$ such that
$\gamma(J_{\lambda}^{-\epsilon}) \geq k$. Moreover, $ 0 \notin
J_{\lambda}^{-\epsilon}$ and $J_{\lambda}^{-\epsilon}\in
\Gamma_{k}$. On the other hand
$$
\displaystyle\sup_{u\in J_{\lambda}^{-\epsilon}}J_{\lambda}(u)\leq
-\epsilon.
$$

Hence,
$$
-\infty < c_{k}=\displaystyle\inf_{C\in
\Gamma_{k}}\displaystyle\sup_{u \in C}J_{\lambda}(u) \leq
\displaystyle\sup_{u\in J_{\lambda}^{-\epsilon}}J_{\lambda}(u) \leq
-\epsilon <0.
$$
\hfill\rule{2mm}{2mm}

The next Lemma allows us to prove the existence of critical points
of $J_{\lambda}$.

\begin{lemma}\label{minimax2}
If $c=c_{k}=c_{k+1}=...=c_{k+r}$ for some $r \in \mathbb{N}$, then
there exists $\lambda^{*}>0$ such that
$$
\gamma(K_{c})\geq r+1,
$$
for $\lambda \in ( 0, \lambda^{*})$.
\end{lemma}
\noindent\textbf{Proof:} Since $c=c_{k}=c_{k+1}=...=c_{k+r} <0$, for
$\lambda^{*}=\min\{\tau_{1},\tau_{2}\}$ and for all $\lambda \in
(0,\lambda^{*})$, from Lemma \ref{nivelbaixo} and Lemma
\ref{minimax1}, we get that $K_{c}$ is compact. Moreover,
$K_{c}= - K_{c}$. If $\gamma(K_{c})\leq r$, there exists a closed
and symmetric set $U$ with $ K_{c}\subset U$ such that $\gamma(U)=
\gamma(K_{c}) \leq r$. Note that we can choose $U\subset
J_{\lambda}^{0}$ because $c<0$. By the deformation lemma
\cite{benci} we have an odd homeomorphism $ \eta: H\rightarrow H$
such that $\eta(J_{\lambda}^{c+\delta}-U)\subset
J_{\lambda}^{c-\delta}$ for some $\delta > 0$ with $0<\delta < -c$.
Thus, $J_{\lambda}^{c+\delta}\subset J_{\lambda}^{0}$ and by
definition of $c=c_{k+r}$, there exists $A \in \Gamma_{k+r}$ such
that $\displaystyle\sup_{u \in A} < c+\delta$, that is, $A \subset
J_{\lambda}^{c+\delta}$ and
\begin{eqnarray}\label{estrela1}
\eta(A-U) \subset \eta ( J_{\lambda}^{c+\delta}-U)\subset
J_{\lambda}^{c-\delta}.
\end{eqnarray}
But $\gamma(\overline{A-U})\geq \gamma(A)-\gamma(U) \geq k$ and
$\gamma(\eta(\overline{A-U}))\geq  \gamma(\overline{A-U})\geq k$.
Then $\eta(\overline{A-U}) \in \Gamma_{k}$ which contradicts
(\ref{estrela1}).\hfill\rule{2mm}{2mm}

\section{Proof of Theorem \ref{Teorema1}}

If $-\infty< c_{1} < c_{2} < ...< c_{k}< ...<0$ with $c_{i}\neq
c_{j}$, once each $c_{k}$ is a critical value of $J_{\lambda}$, we obtain infinitely many critical points of $J_{\lambda}$ and then,
$(T_{\lambda})$ has infinitely many solutions.

On the other hand, if $c_{k}=c_{k+r}$ for some $k$ and $r$, then
$c=c_{k}=c_{k+1}=...=c_{k+r}$ and from Lemma \ref{minimax2}, there
exists $\lambda^{*}>0$ such that
$$
\gamma(K_{c})\geq r+1 \geq 2
$$
for all $\lambda \in (0,\lambda^{*})$. From Proposition
\ref{paracompletar} $K_{c}$ has infinitely many points, that is, $(T_{\lambda})$ has infinitely many solutions.

Let $\lambda^{*}$ be as in Lemma \ref{minimax2} and, for
$\lambda \in (0,\lambda^*)$, let $u_{\lambda}$ be a solution of $(T_{\lambda})$. Thus
$J_{\lambda}(u_{\lambda})=I_{\lambda}(u_{\lambda}) <0$. Hence,
\begin{eqnarray*}
\|u_{\lambda}\|^{2} \leq R_{0},
\end{eqnarray*}
which together with (\ref{comportamentoassimtotico1}) implies that
\begin{eqnarray}\label{C1}
\displaystyle\lim_{\lambda \rightarrow 0}\|u_{\lambda}\|=0.
\end{eqnarray}

Now we use the Moser iteration technique in order to prove that there exists a constant positive $C$, independent on
$\lambda$ such that
\begin{eqnarray}\label{INFTY}
\|u_{\lambda}\|_{\infty} \leq C \|u_{\lambda}\|.
\end{eqnarray}

Using (\ref{INFTY}) we can conclude that
\begin{equation}
\displaystyle\lim_{\lambda\rightarrow 0}\|u_{\lambda}\|_{\infty}=0.
\label{Linfinity}
\end{equation}

In order to save notation, from now on we denote $u_{\lambda}$ by $u$. In
what follows, we fix $R>R_{1}>0$, $R>1$ and take a cut-off function
$\eta_{R} \in C^{\infty}_{0}(\Omega)$ such that $0 \leq \eta_R \leq
1$, $\eta_{R} \equiv  0$ in $B_{R}^{c}$, $\eta_R \equiv 1$ in
$B_{R_{1}}$ and $|\nabla \eta_{R}|\leq C/R$, where $B_{R}\subset
\Omega$ and $C>0$ is a constant.

Let $h(t)= \lambda t^{q-1}+ t^{2^{*}-1}$. Thus
$$
|h(t)|\rightarrow 0 \  \  \mbox{as} \ \  t\rightarrow 0
$$
and
$$
\frac{|h(t)|}{t^{2^{*}-1}} \rightarrow 1 \  \  \mbox{as} \  \
t\rightarrow \infty.
$$

Thus, for all $\delta
>0$ there is $C_{\delta}(\lambda)>0$ such that
\begin{eqnarray}\label{1moser}
h(t)\leq \delta  + C_{\delta}(\lambda) t^{2^{*}-1}.
\end{eqnarray}
Moreover, for $\lambda \in [0,\lambda_0]$, $C_\delta(\lambda)$ can be chosen uniformly in $\lambda$ in such a way that (\ref{1moser}) holds independently of $\lambda$.
For each  $L
> 0$, define
\begin{eqnarray*}
u_{L}(x) =\ \  \left\{
             \begin{array}{l}
              u(x), \quad \mbox{if} \quad u(x) \leq L
        \\
        \\
        L, \quad \mbox{if} \quad u(x) \geq L,
        \\

             \end{array}
           \right.
\end{eqnarray*}
\begin{eqnarray*}
z_{L} = \eta_{R}^{2}u_{L}^{2(\sigma - 1)}u \quad \mbox{and} \quad
w_{L} = \eta_{R}u u_{L}^{\sigma - 1}
\end{eqnarray*}
with $\sigma > 1$ to be determined later. In the course of this
proof, $C_1$, $C_2$..., denote constants independent of
$\lambda$.

\noindent Taking $z_{L}$ as a test function we obtain
$$
I'_{\lambda}(u)z_{L}=0.
$$
More specifically,
$$
\int_{\Omega}a(|\nabla u|^{2})\nabla u \nabla z_{L} = \lambda
\int_{\Omega}u^{q-1}z_{L}+\int_{\Omega}u^{2^{*}-1}z_{L}.
$$
Hence
$$
K_{0}\int_{\Omega}\nabla u \nabla z_{L} \leq \int_{\Omega}h(u)z_{L}.
$$
By \eqref{1moser} we obtain

$$
\int_{\Omega}\nabla u \nabla z_{L} \leq \delta K^{-1}_{0}\int_\Omega z_L +
K^{-1}_{0}C_{\delta}\int_{\Omega}u^{2^{*}-1}z_{L}.
$$
Let us fix $\delta > 0$ small enough in such a way that
$$
\int_{\Omega}\nabla u \nabla z_{L} \leq
C\int_{\Omega}u^{2^{*}-1}z_{L}.
$$
Using $z_{L}$ we obtain
\begin{eqnarray*}
\int_{\Omega}\eta_{R}^{2}u_{L}^{2(\sigma -1)}|\nabla u|^{2} \ dx
&\leq &- \int_{\Omega}\eta_{R}uu_{L}^{2(\sigma-1)}\nabla
\eta_{R}\nabla u \ dx \\
&-& 2(\sigma-1)\int_{\Omega}u_{L}^{(2\sigma-3)}u\nabla u \nabla
u_{L} + \int_{\Omega}\eta_{R}^{2}u^{2^{*}}u_{L}^{2(\sigma-1)}\ dx,
\end{eqnarray*}
and the definition of $u_{L}$ implies
$$
-2(\sigma-1)\int_{\Omega}u_{L}^{(2\sigma-3)}u\nabla u \nabla
u_{L}\leq 0.
$$
Thus
\begin{eqnarray*}
\int_{\Omega}\eta_{R}^{2}u_{L}^{2(\sigma -1)}|\nabla u|^{2} \ dx\leq
+ \int_{\Omega}\eta_{R}uu_{L}^{2(\sigma-1)}|\nabla \eta_{R}| |\nabla
u| \ dx  + \int_{\Omega}\eta_{R}^{2}u^{2^{*}}u_{L}^{2(\sigma-1)}\
dx.
\end{eqnarray*}

\noindent Taking $z_{L}$ as a test function and using
(\ref{1moser}), we obtain
\begin{eqnarray*}
\int_{\Omega}\eta_{R}^{2}u_{L}^{2(\sigma -1)}|\nabla u|^{2} \ dx\leq
C_{1} \int_{\Omega}\eta_{R}uu_{L}^{2(\sigma-1)}|\nabla
\eta_{R}||\nabla u| \ dx + C_{1}
\int_{\Omega}\eta_{R}^{2}u^{2^{*}}u_{L}^{2(\sigma-1)}\ dx.
\end{eqnarray*}
Fixing $\widetilde{\tau}>0$ and using Young's inequality, we obtain
\begin{eqnarray*}
\int_{\Omega}\eta_{R}^{2}u_{L}^{2(\sigma -1)}|\nabla u|^{2} \ dx
\leq && C_{1} \int_{\Omega}\bigg(\widetilde{\tau}\eta_{R}^{2}|\nabla
u|^{2}+ C_{\widetilde{\tau}}u^{2}|\nabla \eta_{R}|^{2}
\biggl)u_{L}^{2(\sigma-1)} \ dx + \\ && C_{1}
\int_{\Omega}\eta_{R}^{2}u^{2^{*}}u_{L}^{2(\sigma-1)}\ dx.
\end{eqnarray*}
Choosing $\widetilde{\tau} \leq 1/4$, it follows that
\begin{eqnarray}\label{2moser}
\int_{\Omega}\eta_{R}^{2}u_{L}^{2(\sigma -1)}|\nabla u|^{2} \ dx\leq
C_{2}\bigg(\int_{\Omega}u^{2}u_{L}^{2(\sigma-1)}|\nabla
\eta_{R}|^{2}  \ dx +
\int_{\Omega}\eta_{R}^{2}u^{2^{*}}u_{L}^{2(\sigma-1)}\ dx\bigg).
\end{eqnarray}

\noindent On the other hand, we get
\begin{eqnarray*}
S \|w_{L}\|^{2}_{L^{2^{*}}(\Omega)} & \leq &
 \int_{\Omega} \left|\nabla
\left(\eta_{R}u u_{L}^{\sigma-1}\right)\right|^2  \\
& \leq &   \int_{\Omega} |u|^2 u_{L}^{2(\sigma-1)}|\nabla
\eta_{R}|^2 +
 \int_{\Omega} \eta_{R}^2\left|\nabla \left(u
u_{L}^{\sigma-1}\right)\right|^2.
\end{eqnarray*}
But
\begin{eqnarray*}
 \int_{\Omega} \eta_{R}^2\left|\nabla \left(u
u_{L}^{\sigma-1}\right)\right|^2 &=&
 \int_{\{|u| \leq L\}}\eta_{R}^2\left|\nabla
\left(u u_{L}^{\sigma-1}\right)\right|^2  +
 \int_{\{|u|
> L\}} \eta_{R}^2\left|\nabla
\left(u u_{L}^{\sigma-1}\right)\right|^2 \\
&=&  \int_{\{|u| \leq L\}} \eta_{R}^2\left|\nabla
u^{\sigma}\right|^2 +
 \int_{\{|u| > L\}} \eta_{R}^2
L^{2(\sigma-1)}\left|\nabla u\right|^2
\\
& \leq & \sigma^2  \int_{\Omega}\eta_{R}^2
u_{L}^{2(\sigma-1)}|\nabla u|^2,
\end{eqnarray*}
and therefore
$$
\|w_{L}\|^{2}_{L^{2^{*}}(\Omega)} \leq C_3\sigma^{2} \left(
\int_{\Omega}
  |u|^{2}u_{L}^{2(\sigma - 1)}|\nabla \eta_{R}|^{2} +
 \int_{\Omega}  \eta_{R}^{2} u_{L}^{2(\sigma - 1)}|\nabla
u|^{2}\right).
$$
From this and (\ref{2moser}),
\begin{equation}\label{4moser}
\|w_{L}\|^{2}_{L^{2^{*}}(\Omega)} \leq C_4\sigma^{2} \left(
\int_{\Omega} |u|^{2}u_{L}^{2(\sigma - 1)}|\nabla \eta_{R}|^{2} +
\int_{\Omega} \eta_{R}^{2}|u|^{2^{*}}u_{L}^{2(\sigma - 1)}\right),
\end{equation}
for all $\sigma>1$. The above expression, the properties of
$\eta_{R}$ and $u_{L} \leq u$, imply that
\begin{equation} \label{5moser}
\|w_{L}\|^{2}_{L^{2^{*}}(\Omega)} \leq C_4\sigma^{2} \int_{B_{R}}
\left( |u|^{2\sigma}|\nabla \eta_{R}|^{2} +
|u|^{2^{*}-2}|u|^{2\sigma}\right).
\end{equation}

If we set
\begin{equation} \label{def_t}
t:= \frac{2^*2^*}{2(2^*-2)}>1,~~~~\alpha := \frac{2t}{t-1}< 2^*,
\end{equation}
we can apply H\"{o}lder's inequality with exponents $t/(t-1)$ and
$t$ in (\ref{5moser}) to get
\begin{equation} \label{6moser}
\begin{array}{lcl}
\|w_{L}\|^{2}_{L^{2^{*}}(\Omega)} &\leq& C_4\sigma^{2}
\|u\|_{L^{\sigma\alpha}(B_{R})}^{2\sigma} \left( \int_{B_{R}}|\nabla
\eta_{R}|^{2t}\right)^{1/t}
 \\
&& +C_4\sigma^{2} \|u\|_{L^{\sigma\alpha}(B_{R})}^{2\sigma} \left(
\int_{B_{R}}|u|^{2^*(2^*/2)}\right)^{1/t}.
\end{array}
\end{equation}
Since $\eta_{R}$ is constant on $B_{R_1} \cup B_{R}^c$ and $|\nabla
\eta_{R}| \leq C/R$, we conclude that
\begin{equation} \label{61moser}
 \int_{B_{R}}|\nabla \eta_{R}|^{2t} =
 \int_{B_{R}\backslash B_{R_{1}}}|\nabla \eta_{R}|^{2t}
\leq \frac{C_5}{R^{2t-N}} \leq C_5.
\end{equation}
We have used $R>1$ and $2t=\frac{2^{*}}{2}N > N$ in the last
inequality.

\noindent \textbf{Claim.} There exist  a constants $K >0$
independent on $\lambda$ such that,
$$
\int_{\Omega}|u|^{2^*(2^*/2)}\leq K.
$$

Assuming the claim is true, we can use (\ref{6moser}) and
(\ref{61moser}) to conclude that
\begin{eqnarray*}
\|w_{L}\|^{2}_{L^{2^{*}}(\Omega)} \leq
C_6\sigma^{2}\|u\|^{2\sigma}_{L^{\sigma\alpha}(B_{R})}.
\end{eqnarray*}
Since
\begin{eqnarray*}
\|u_{L}\|^{2\sigma}_{L^{\sigma 2^{*}}(B_R)} &=& \left(
\int_{B_R}u_{L}^{\sigma 2^{*}
}\right)^{2/2^*} \\
&\leq&   \left(  \int_{\Omega} \eta_{R}^{2^*}
|u|^{2^{*}}u_{L}^{2^{*}(\sigma -1) }\right)^{2/2^*}  \\
&=&\|w_{L}\|^{2}_{L^{2^{*}}(\Omega)} \leq
C_6\sigma^{2}\|u\|^{2\sigma}_{L^{\sigma\alpha}(\Omega)},
\end{eqnarray*}
we can apply Fatou's lemma in the variable $L$ to obtain
\begin{eqnarray*}
\|u\|_{L^{\sigma 2^{*}}(B_R)}\leq
C_7^{1/\sigma}\sigma^{1/\sigma}\|u\|_{L^{\sigma\alpha}(\Omega)},
\end{eqnarray*}
whenever $|u|^{\sigma\alpha} \in L^1(B_{R})$. Here, $C_7$ is a
positive constant independent on $R$. Iterating this process, for
each $k\in\mathbb{N}$, it follows that
$$
\|u\|_{L^{\sigma^k 2^{*}}(B_R)} \leq C_7^{\sum_{i=1}^k \sigma^{-i}}
\sigma^{\sum_{i=1}^m i \sigma^{-i}} \|u\|_{L^{2^*}(\Omega)}.
$$
Since $\Omega$ can be covered by a finite number of balls $B_{R}^j$,
we have that
$$
\|u\|_{L^{\sigma^k 2^{*}}(\Omega)}\leq \displaystyle\sum_j^{finite}
\|u\|_{L^{\sigma^k 2^{*}}(B_R^j)} \leq \displaystyle\sum_j^{finite}
C_7^{\sum_{i=1}^k \sigma^{-i}} \sigma^{\sum_{i=1}^m i \sigma^{-i}}
\|u\|_{L^{2^*}(\Omega)} .
$$
Since $\sigma>1$,  we let $k \to \infty$ to get
$$
\|u\|_{L^{\infty}(\Omega)} \leq K_{2}\|u\|,
$$
for some $K_{2}>0$ independent on $\lambda$.

It remains to prove the claim. From (\ref{4moser})
\begin{equation}\label{7moser}
\|w_{L}\|^{2}_{L^{2^{*}}(\Omega)} \leq C_{9}\sigma^{2} \left(
\int_{\Omega} |u|^{2}u_{L}^{2(\sigma - 1)}|\nabla \eta_{R}|^{2} +
\int_{\Omega} \eta_{R}^{2}|u|^{2^{*}}u_{L}^{2(\sigma - 1)}\right),
\end{equation}
We set $\sigma:= 2^*/2$ in \eqref{4moser} to obtain
$$
\|w_{L}\|^{2}_{L^{2^{*}}(\Omega)} \leq C_{10}\left(
 \int_{\Omega}
|u|^{2}u_{L}^{(2^{*}-2)}|\nabla \eta_{R}|^{2} + \int_{B_{R}}
\eta_{R}^{2}|u|^{2}u_{L}^{(2^{*}-2)}|u|^{(2^{*}-2)}\right).
$$
By H\"{o}lder's inequality with exponents $2^*/2$ and $2^*/(2^*-2)$
we get
\begin{eqnarray*}
\|w_{L}\|^{2}_{L^{2^{*}}(\Omega)} &\leq & C_{10} \int_{\Omega}
|u|^{2}u_{L}^{(2^{*}-2)}|\nabla \eta_{R}|^{2}\\&+& C_{10} \left(
\int_{B_{R}}\left(\eta_{R}|u|
u_{L}^{(2^*-2)/2}\right)^{2^{*}}\right)^{2/2^*}
\|u\|_{L^{2^*}(\Omega)}^{2^*-2}.
\end{eqnarray*}
From (\ref{C1}) and  recalling that $\eta_{R} u
u_{L}^{(2^*-2)/2}=w_{L}$, $u_{L} \leq u$ and $\nabla \eta_R$ is
bounded, we obtain
\begin{eqnarray*}
\|w_{L}\|^{2}_{L^{2^{*}}(\Omega)} \leq C_{11} \int_{\Omega}
|u|^{2}u_{L}^{(2^{*}-2)}|\nabla \eta_{R}|^{2} \leq C_{11}
\int_{\Omega} |u|^{2^*} \leq C_{12}.
\end{eqnarray*}
The definition of $\eta_{R}$ and $w_{L}$ and the above inequality
imply that
$$
\left( \int_{B_{R}} |u|^{2^*}u_{L}^{2^*(2^*-2)/2}\right)^{2/2^*}
\leq |w_{L}|^{2}_{L^{2^{*}}(\Omega)} \leq C_{12}.
$$
Using Fatou's lemma in the variable $L$, we have
\begin{eqnarray*}
\int_{B_{R}}|u|^{2^*(2^*/2)}\leq K:=C_{12}^{2^*/2}.
\end{eqnarray*}
Since $\Omega$ can be covered by a finite number of balls $B^j_{R}$,
we have that
$$
\int_{\Omega}|u|^{2^*(2^*/2)}\leq
\displaystyle\sum_j^{finite}\int_{B_{R}}|u|^{2^*(2^*/2)}\leq K_{3},
$$
for some $K_{3}>0$.

In order to estimate $\|\nabla u_\lambda\|_ \infty$, we make use of the following result by Stampacchia in \cite{Stampacchia}
\begin{lemma}
Let $A(\eta)$ a given $C^1$ vector field in $\mathbb{R}^N$, and $f(x,s)$ a bounded Carath\'eodory function in $\Omega \times \mathbb{R}$. Let $u \in H^1_0(\Omega)$ be a solution of
$$\int_\Omega \left(A(|\nabla u|)\nabla \varphi + f(x,u)\varphi\right) = 0,$$
for all $\varphi \in H^1_0(\Omega)$. Assume that there exist $0 < \nu < M$ such that
\begin{equation}
\nu |\xi|^2 \leq \frac{\partial A_i}{\partial \eta_j}(\nabla u)\xi_i\xi_j, \quad \mbox{and} \quad \left|\frac{\partial A_i}{\partial \eta_j}(\nabla u)\right| \leq M,
\label{elipticidade}
\end{equation}
for all $i,j = 1, ... , N$ and $\xi \in \mathbb{R}^N$. Then $u \in W^{2,p}(\Omega) \cap C^{1,\alpha}(\overline{\Omega})$, for all $\alpha \in (0,1)$ and $p > 1$. Moreover
\begin{equation}
\|u\|_{1,\alpha} \leq C(\nu, M, \Omega) \|f(\cdot,u)\|_\infty.
\end{equation}
\end{lemma}

By the definition of $a$, for $r$ small enough (\ref{elipticidade}) hold. This, together with the fact that $\|u_\lambda\|_\infty$ is bounded allow us to apply the last result. Then (\ref{Linfinity}) implies that
\begin{equation}
\|u\|_{1,\alpha} \leq \lambda \|u\|_\infty^{q-1} + \|u\|_\infty^{2^*-1} = o(\lambda),
\label{limitegradiente}
\end{equation}
as $\lambda \to 0$.

Then, there exists $\lambda^* > 0$ such that $\lambda \in (0,\lambda^*)$ implies that $\|\nabla u\|_\infty \leq r$ and hence, $u_\lambda$ is a solution of $(P_\lambda)$.
\hfill\rule{2mm}{2mm}

\end{document}